\documentclass[10pt]{amsart}

\usepackage{latexsym,amssymb,amsmath,graphics}%,showkeys}
\usepackage[dvips]{graphicx}

\def\ff{{\mathcal F}}

\def\bpf{\begin{proof}[Proof]}
\def\epf{\end{proof}}
\def\ffi{\varphi}

\def\dst{\displaystyle}

\def\C{{\mathbb{C}}}

\def\R{{\mathbb{R}}}

\newcommand{\norm}[1]{{\left\|{#1}\right\|}}
\newcommand{\ent}[1]{{\left[{#1}\right]}}
\newcommand{\abs}[1]{{\left|{#1}\right|}}
\newcommand{\scal}[1]{{\left\langle{#1}\right\rangle}}
\renewcommand{\iint}{\int\!\!\!\int}
\newenvironment{remark}[1][]{\vskip1pt\noindent\rm\textit{Remark}\,:\ }{\rm\vskip1pt}

\newtheorem{lem}{Lemma}[section]
\newtheorem{prop}[lem]{Proposition}
\newtheorem{theo}[lem]{Theorem}
\newtheorem{cor}[lem]{Corollary}

\newenvironment{lemma}[1][]{\lem[#1]\sl}{\rm\vskip3pt}
\newenvironment{proposition}[1][]{\prop[#1]\sl}{\rm\vskip3pt}
\newenvironment{theorem}[1][]{\theo[#1]\sl}{\rm\vskip3pt}
\newenvironment{corollary}[1][]{\cor[#1]\sl}{\rm\vskip3pt}

\begin{document}

\title[uncertainty principles for the Fourier and the windowed Fourier
transforms]{Hermite functions and uncertainty
principles for the Fourier and the windowed Fourier transforms}
\author{Aline Bonami, Bruno Demange \& Philippe Jaming}
\address{Universit\'e d'Orl\'eans\\
Facult\'e des Sciences\\
D\'epartement de Math\'ematiques\\BP 6759\\ F 45067 ORLEANS Cedex 2\\
FRANCE}
\email{bonami@labomath.univ-orleans.fr\\
jaming@labomath.univ-orleans.fr}

\address{E.N.S. Lyon\\ 46 All\'ee d'Italie\\ 60364 Lyon CEDEX 07\\
FRANCE\\}
\email{bdemange@ens-lyon.fr}

\begin{abstract}
We  extend an uncertainty principle due to  Beurling into a
characterization of Hermite functions. More precisely,   all
functions $f$ on $\R^d$ which may be written as $P(x)\exp (Ax,
x)$, with $A$ a real symmetric definite positive matrix, are
characterized by integrability conditions on the product $f(x)
\widehat { f}(y)$. We also give the best constant in uncertainty
principles of Gelf'and Shilov type.
 We then  obtain similar results
for the windowed Fourier transform (also known, up to elementary
changes of functions,  as the radar ambiguity function
or the Wigner transform). We complete the paper with a  sharp version of
Heisenberg's inequality for this transform.
\end{abstract}

\keywords{Uncertainty principles; short-time Fourier transform;
windowed Fourier transform; Gabor transform; ambiguity function;
Wigner transform; spectrogramm}
\subjclass{42B10;32A15;94A12}

\thanks{Research partially financed by : {\it European Commission}
(TMR 1998-2001 Network {\it Harmonic Analysis}).}
\maketitle

\section{Introduction and Notations.}

Uncertainty principles state that a function and its Fourier transform
cannot be simultaneously sharply localized. To be more precise, let
$d\geq 1$ be the dimension, and let us denote by $\scal{.,.}$ the scalar product
 and by $\|.\|$ the Euclidean norm on $\R^d$. Then, for $f\in
L^2(\R^d)$, define the Fourier transform of $f$ by
$$
\widehat f(y)=\int_{\R^d}f(t)e^{-2i\pi\scal{t,y}}dt.
$$
The most famous uncertainty principle, due to Heisenberg and Weil, can be stated
in the following directional version\,:

{\bf Heisenberg's inequality.}
%\label{th2}
Let $i=1,\ldots,d$ and $f\in L^2(\R^d)$. Then
\begin{equation}
\inf_{a\in\R}\left ({\int_{\R^d}{(x_i-a)^2\abs{f(x)}^2dx}}\right )
\inf_{b\in\R}\left
({\int_{\R^d}{(\xi_i-b)^2\abs{\widehat{f}(\xi)}^2d\xi}}
\right )\geq\frac{{\|f\|} _{L^2} ^4}{16\pi^2}.
\label{eqth2}
\end{equation}
Moreover (\ref{eqth2}) is an equality if and only if $f$ is of the
form $$ f(x)=C(x_1,\dots,x_{i-1},x_{i+1},\dots,x_n)e^{-2i\pi
bx_i}e^{-\alpha(x_i-a)^2} $$ where  $C$ is a function in
$L^2(\R^{d-1})$, $\alpha>0$, and $a$ and $b$ are real constants
for which  the two infimums in (\ref{eqth2}) are realized.

The usual non-directional uncertainty principle follows easily
from this one.  We refer to  the recent survey articles by Folland
and Sitaram \cite{FS} and Dembo, Cover and Thomas \cite{DCT} as
well as the book of Havin and J\"oricke \cite{HJ} for various
uncertainty principles of different nature which may be found in
the literature. One theorem stated in \cite{FS} is due to
Beurling. Its proof has been written much later by H\"ormander in
\cite {HB}.  Our first aim  is to weaken the assumptions so that
 non zero solutions given by Hermite
functions are also possible.
 More precisely, we will
prove the following theorem\,:

\begin{theorem}[Beurling-H\"ormander type] Let $f\in L^2(\R^d)$ and $N\geq 0$.  Then
\begin{equation}
\iint_{\R^d\times\R^d}
\frac{|f(x)||\widehat{f}(y)|}{(1+\norm{x}+\norm{y})^N}
e^{2\pi\abs{\scal{x,y}}}dxdy<+\infty
\label{condBHT}
\end{equation}
if and only if $f$ may be written as
$$
f(x)=P(x)e^{-\pi\scal {Ax, x}},
$$
where $A$ is a real
positive definite symmetric matrix and $P$ is a polynomial of degree
$<\frac{N-d}{2}$.
\label{BHT}
\end{theorem}
In particular, for $N\leq d$, the function $f$ is identically $0$.
 Beurling-H\"ormander's original theorem is the above theorem for
$d=1$ and $N=0$. An extension to $d\geq 1$ but still $N=0$ has
been given, first by S.C. Bagchi and S. K. Ray  in \cite{BaSw} in
a weaker form, then very recently by S. K. Ray and E. Naranayan in
the present form. Their proof, which relies on the one dimensional
case, uses Radon transform \cite{RN}.

Let us remark that it is
easy to adapt  the usual proof of
  Hardy's Theorem to allow Hermite
functions (see for instance  \cite {Pfann}), while the proof given
here for integrability conditions uses new ingredients. At the same
time, it simplifies H\"ormander's argument for the case
$N=0$, $d=1$, in such a way
that the proof can now be given in any textbook on Fourier Analysis. We
give the simplified proof in the Appendix, since it may be useful in
this context.

The previous theorem has as an immediate corollary the following
characterization.
\begin{corollary}
A function $f\in L^2(\R^d)$ may be written as $$
f(x)=P(x)e^{-\pi\scal {Ax, x}}, $$ with $A$  a real positive definite
symmetric matrix and $P$  a polynomial, if and only if the
function $$f(x)\widehat {f}(y)\exp(2\pi\abs{\scal{x,y}})$$ is
slowly increasing on $\R^d\times\R^d$. \label{HBP3}\end{corollary}
As an easy consequence of the previous theorem, we also deduce the
following corollary, which generalizes the Cowling-Price
uncertainty principle (see \cite {CP}).
\begin{theorem}[Cowling-Price type]
Let $N\geq 0$. Assume that $f\in L^2(\R^d)$ satisfies
$$
    \int_{\R^d}|f(x)|\frac{e^{\pi a\abs{x_{j}}^2}}{(1+\abs{x_{j}})^N}dx<+\infty
\hspace{2ex}\mathrm{and} \hspace{2ex}
\int_{\R^d}|\widehat{f}(y)|\frac{e^{\pi
b\abs{y_{j}}^2}}{(1+\abs{y_{j}})^N} dy<+\infty $$ for $j=1,
\cdots, d$ and for some positive constants $a$ and $b$ with
$ab=1$. Then $f(x)=P(x)e^{-a\norm{x}^2}$ for some polynomial $P$.
\label{Cowl1}\end{theorem}

The Cowling-Price theorem is given in  \cite {CP} for $d=1$,
$N=0$, and with $p$-th powers as well. This last extension is
a trivial consequence of H\"older's inequality once $N$ is allowed
to take positive values.

It is remarked in  \cite {HB}, as well as in \cite{BaSw}, that a
theorem of Beurling-H\"ormander type implies also a theorem of Gel'fand-Shilov
type. In fact, the constant that one obtains when doing this is not the best one, and
looking carefully at the literature on entire functions one gets the
following theorem, which gives the critical constant.
\begin{theorem}[Gel'fand-Shilov type]
Let $1<p<2$, and let $q$ be the conjugate exponent. Assume that
$f\in L^2(\R^d)$ satisfies $$
\int_{\R^d}|f(x)|e^{2\pi\frac{a^p}{p}\abs{x_{j}}^p}dx<+\infty
\hspace{2ex} \mathrm{and} \hspace{2ex}
\int_{\R^d}|\widehat{f}(y)|e^{2\pi\frac{b^q}{q}\abs{y_{j}}^q}dy<+\infty
$$ for some $j=1, \cdots, d$ and for some positive constants $a$
and $b$. Then $f=0$ if $ab>\abs {\cos(\frac{p\pi}{2})}^{\frac
1p}$. If $ab< \abs {\cos(\frac{p\pi}{2})}^{\frac 1p}$, one may
find a dense subset of  functions which satisfy the above
conditions for all $j$. \label{GelShi}
\end{theorem}
From Theorem \ref{BHT}, one only gets the first result for
$ab\geq1$, which is clearly much weaker. Close but different
statements are given in \cite{HJ}.
\medskip

One way one may hope to overcome the lack of localisation is to
use the {\it windowed Fourier transform}, also known as the
(continuous) {\it Gabor transform} or the {\it short-time Fourier
transform}. To be more precise, fix $v\in L^2(\R^d)$, the
``window'', and define for $u\in L^2(\R^d)$\,: $$
S_vu(x,y)=\overline{\widehat{\overline{u}v(.-x)}}(y)
=\int_{\R^d}u(t)\overline{v(t-x)}e^{2i\pi\scal{t,y}}dt. $$ This
transform occurs also in several other forms, for example
$\abs{S_vu}^2$ is known as a {\it spectrogram}. For sake of
symmetry in $u$ and $v$, we rather focus on the {\it radar
ambiguity function} defined for $u,v\in L^2(\R^d)$ by
\begin{equation}
A(u,v)(x,y)=\int_{\R^d}u\left(t+\frac{x}{2}\right)
\overline{v\left(t-\frac{x}{2}\right)}e^{-2i\pi\scal {y, t}}dt.
\label{defafg}
\end{equation}
Since $\abs{A(u,v)}=\abs{S_v u}$, there will be no loss in doing
so. We refer the reader to \cite{AT}, \cite{CB}, \cite{Co} and the
references there for the way these functions occur in signal
processing, and their basic properties.

Finally, $W(u,v)$, the Fourier transform of $A(u,v)$ in $\R^{2d}$ is
known in quantum mechanics and in the PDE community as the {\it Wigner
transform} or {\it Wigner distribution}. Since
$$
W(u,v)(x,y)=\int_{\R^d}u\left(x+\frac{t}{2}\right)
\overline{v\left(x-\frac{t}{2}\right)}e^{2i\pi \scal {y, t}}dt,
$$
$W(u,v)$ is also related to $A(u,v)$ by
$$
W(u,v)(x,y)=2^dA(u,Zv)(2x,-2y)
$$
where $Zv(x)=v(-x)$. So again, all results stated here can be restated
in terms of the Wigner transform.

Our second aim here is to extend  uncertainty principles to the
radar ambiguity functions. In particular, we will show that
$A(u,v)$ satisfies theorems of Cowling-Price type on one side, of
Gel'fand-Shilov type on the other one. Both results are sharp,
with the same characterization of the Hermite functions in the
first case.

We also give a version of Heisenberg's inequality for $A(u,v)$ that
is stronger than a previous version
by A.J.E.M. Janssen (\cite{janss1}, see also \cite{flandrin}). The one
dimensional case for the Wigner transform $W(u,u)$ can be found in
\cite{dBr}. This Heisenberg's inequality may
be stated in the following matricial form.
\begin{theorem}\label{variance}
Assume that $u$ and $v$ be in $L^2(\R^d)$, with
$\norm{u}_{L^2}\norm{v}_{L^2}1$ and $$ \int_{R^d}\norm{x}^2
(\abs{u}^2+\abs{v}^2+\abs{\widehat u}^2+\abs{\widehat
v}^2)dx<\infty\,. $$
    Let $(X,Y)$ be a random vector  with probability density given on
    $\R^d\times\R^d$ by the function $|A(u,v)|^2$. Then $X$ and $Y$
    are not correlated, and their two
    covariance matrices are such that
    $$4\pi^2 V(X)-V(Y)^{-1}$$
    is semi-positive definite. Moreover, if it is the zero matrix,
    then $u$ and $v$ are Gaussian functions.
    \end{theorem}
    In Radar Theory, for which $d=1$, the couple $(X,Y)$ may be given a
    physical meaning: its first component is related to the distance
    to the target, the second to its velocity. So the variance gives
    the quadratic error when estimating the distance or the velocity
    by the corresponding mean.

A different problem consists in minimizing the same quantity
$A(u,v)$ for a fixed ``window" $v$, or more generally to consider
uncertainty principles in terms of a function $u$ and its windowed
Fourier transform $A(u,v)$. Such problems have been considered by
E. Wilczock \cite{Wilc}. Such results follow from ours only when
the window is Gaussian.

The article is organized as follows\,: the next section is devoted
to the proof of Theorem \ref{BHT}, whereas in section
\ref{corollaires}, we consider the other mentioned uncertainty
principles for the Fourier transform. In section \ref{propamb}, we
recall some basic properties of the ambiguity functions, pursuing
in section \ref{heisamb} with the Heisenberg inequality for these
functions. Section \ref{BHamb} is devoted to the extension of
uncertainty principles to ambiguity functions.

\section{Generalization of Beurling-H\"ormander's theorem.}
The statement of Theorem \ref{BHT} is divided into two
propositions. The first one is elementary.

\begin{proposition} Let $f\in L^2(\R^d)$ be a function of the form
$$
f(x)=P(x)e^{-\pi \scal {(A+iB)x, x}},
$$
with  $A$ and $B$ two real
symmetric matrices and $P$  a polynomial. Then the matrix $A$ is
positive definite. Moreover, the
three following conditions are equivalent:
\begin{itemize}
\item [$(i)$]\hspace{0.5cm}
$B=0$ and $\deg (P)<\frac{N-d}{2}$\ ;
\item [$(ii)$]\hspace{0.5cm}
$\dst
\iint_{\R^d\times\R^d}
\frac{|f(x)||\widehat{f}(y)|}{(1+\norm{x}+\norm{y})^{N}}
e^{2\pi\abs{\scal{x,y}}}dxdy<+\infty\ ;$
\item [$(iii)$]\hspace{0.5cm}
$\dst
\iint_{\R^d\times\R^d}
\frac{|f(x)||\widehat{f}(y)|}{(1+\norm{x})^\frac N2(1+\norm{y})^\frac N2}
e^{2\pi\abs{\scal{x,y}}}dxdy<+\infty\ .$
\end{itemize}
\end{proposition}

\bpf
The fact that $A$ is positive definite is elementary. Then,
    after a change of
variables, we may assume that $A$ is  the identity matrix,
so that $f$ may be written as
$P(x) e^{-\pi\norm { x}^2} e^{-i\pi\scal {Bx, x}}$.
The Fourier transform of $f$ may be written as
$Q(y)  e^{-\pi\scal {(I+iB)^{-1}y, y}}$, with $\deg(P)=\deg(Q)=n$.

To prove that $(i)$
implies $(iii)$, it is sufficient to prove that, for $\alpha>\frac
d2$, we have the inequality
$$
\iint_{\R^d\times\R^d}{(1+\norm{x})^{-\alpha}(1+\norm{y})^{-\alpha}}
e^{-\pi\norm{x-y}^2}dxdy<+\infty\ .
$$
But this integral is twice the integral on the subset where
$\norm {x}\leq \norm {y}$. So it is bounded by
$$
2\int_{\R^d}\ent{(1+\norm{x})^{-2\alpha}\times
\int_{\R^d}e^{-\pi\norm{x-y}^2}dy}dx<+\infty
$$
which allows to conclude.

It is clear  that $(iii)$ implies $(ii)$ for all functions $f$.
Let us now prove that $(ii)$ implies $(i)$. First, writing
$(1+iB)^{-1}=(I-iB)(I+B^2)^{-1}$, it is immediate that the
integrability of $f(x)\widehat f(y)
\frac{\exp2\pi\abs{\scal{x,y}}}{(1+\norm{x}+\norm{y})^N}$ implies
that the homogeneous polynomial $\norm{x}^2 +\scal {(I+B^2)^{-1}y,
y}-2\scal{x,y}$ is always positive, which implies that $B=0$. Now,
let $P_{n}$ and $Q_{n}$ be the homogeneous terms of maximal degree
of the polynomials $P$ and $Q$. There exists $x^{(0)}\in \R^d$ of
norm $1$ such that $P_{n}(x^{(0)}) Q_{n}(x^{(0)})$ is different
from $0$. We call $\Gamma_{\varepsilon}$ the cone which is
obtained as the image of the cone $$ \left\{x=(x_{1}, \cdots,
x_{d})\ ;\ (x_{2}^2+\cdots+x_{d}^2)^{\frac{1}{2}}<\varepsilon
x_{1}\right\} $$ under a rotation which maps the point $(1, 0,
\cdots, 0)$ to $x^{(0)}$. Then, for $\varepsilon<1$ small enough,
there exists a constant $c$ such that, for $x, y \in
\Gamma_{\varepsilon}$, one has the inequality $$ |P_{n}(x)|
|Q_{n}(y)|\geq c \norm {x}^n \norm {y}^n\,. $$ The same inequality
is valid for $P$ and $Q$ for $x$ and $y$ large, which implies that
$$ \iint_{\Gamma_{\varepsilon}\times \Gamma_{\varepsilon}}
\frac{\norm {x}^n \norm {y}^n}{(1+\norm{x}+\norm{y})^{N}}
e^{-\pi\norm {x-y}^2}dxdy<+\infty\,. $$ We remark that if
$x\in\Gamma_{\varepsilon}$ then $-\Gamma_{\varepsilon}\subset
x-\Gamma_{\varepsilon}$. So, {\it a fortiori}, we have that $$
\iint_{\Gamma_{\varepsilon}\times (-\Gamma_{\varepsilon})}
\frac{\norm {x}^n \norm {y}^n}{(1+\norm{x}+\norm{t+x})^{N}}
e^{-\pi\norm {t}^2}dxdt<+\infty\,. $$
We know, using Fubini's
theorem, that there exists $t$ for which $$
\int_{\Gamma_{\varepsilon}} \frac{\norm {x}^n \norm
{y}^n}{(1+\norm{x}+\norm{t+x})^{N}} dx <+\infty\,, $$ which proves
that $N-2n>d$. \epf

We have written here the equivalence between Conditions $(ii)$ and
$(iii)$ to remark that Condition $(iii)$ could be written in place
of (\ref{condBHT}) in Theorem \ref {BHT}.

\medskip

The next proposition is much deeper.

\begin{proposition}
 Let $f\in L^2(\R^d)$. If, for some positive integer $N$,
\begin{equation}
\iint_{\R^d\times\R^d}
\frac{|f(x)||\widehat{f}(y)|}{(1+\norm{x}+\norm{y})^N}
e^{2\pi\abs{\scal{x,y}}}dxdy<+\infty\, ,
\label{condBHT2}
\end{equation}
 then $f$ may be written as
$$
f(x)=P(x)e^{-\pi \scal {(A+iB)x, x}},
$$
where $A$ and $B$ are two  symmetric matrices and $P$ is a polynomial.
\label{BHT2}
\end{proposition}

\bpf
We may  assume that $f\not=0$.
\vskip3pt
\noindent{\it First step.} {\sl Both $f$ and $\widehat{f}$
are in $L^1(\R^d)$.}
\vskip3pt
\noindent For almost every $y$,
$$
|\widehat{f}(y)|\int_{\R^d}\frac{|f(x)|}{(1+\norm{x})^N}e^{2\pi
\abs{\scal{x,y}}}dx<+\infty.
$$
As $f\neq 0$, the set of all $y$'s such that $\widehat{f}(y)\neq 0$
has positive measure. In particular, there is a basis
$y^1\dots y^d$ of $\R^d$ such that, for $i=1,\ldots,d$,
$\widehat{f}(y^{(i)})\neq 0$ and
$$
\int_{\R^d}\frac{|f(x)|}{(1+\norm{x})^N}e^{2\pi
\abs{\scal{x,y^{(i)}}}}dx<+\infty.
$$
Since, clearly, there exists a constant $C$ such that
$$
(1+\norm{x})^N\leq C\sum_{i=1}^d \exp\ent{2\pi\abs{\scal{x,y^{(i)}}}},
$$
we get
$f \in L^1(\R^d)$. Exchanging the roles of $f$ and $\widehat{f}$, we get
$\widehat{f}\in L^1(\R^d)$. \hfill$\diamond$

\vskip3pt
\noindent{\it Second step.} {\sl The function $g$ defined by
$\widehat{g}(y)=\widehat{f}(y)e^{-\pi\norm{y}^2}$
satisfies the following properties (with $C$ depending only on $f$)
\begin{equation}
    \int_{\R^d}|\widehat{g}(y)| e^{\pi \norm{y}^2}dy < \infty\,;
\label{step2p0} \end{equation}
\begin{equation}
|\widehat{g}(y)|\leq Ce^{-\pi \norm{y}^2}\,;
\label{step2p1} \end{equation}
\begin{equation}\label{step2p2} \dst\iint_{\R^d\times\R^d}
\frac{|g(x)||\widehat{g}(y)|}{(1+\norm{x}+\norm{y})^N}
e^{2\pi\abs{\scal{x,y}}}dxdy<+\infty\,;
 \end{equation}
 \begin{equation}\label{step2p3} \dst\int_{\norm {x}\leq R}\int_{\R^d}
|g(x)||\widehat{g}(y)|
e^{2\pi\abs{\scal{x,y}}}dxdy<C(1+R)^N\,.
 \end{equation}
} \vskip3pt \noindent Property (\ref{step2p0}) is obvious from the
definition of $g$ and the fact that $\widehat f$ is in
$L^1(\R^d)$. As $f\in L^1(\R^d)$, $\widehat{f}$ is bounded thus
(\ref{step2p1}) is also obvious. To prove (\ref{step2p2}), we have
$$ \iint_{\R^d\times\R^d}
\frac{|g(x)||\widehat{g}(y)|}{(1+\norm{x}+\norm{y})^N}
e^{2\pi\abs{\scal{x,y}}}dxdy \leq\iint
|f(t)||\widehat{f}(y)|A(t,y)e^{2\pi \abs{\scal{t,y}}}dtdy $$ with
$$ A(t,y)=\int \frac{e^{-\pi \norm{x}^2-\pi \norm{y}^2+2\pi
\abs{\scal{x,y}}}}{(1+\norm{t-x}+\norm{y})^N}dx. $$ We claim that
\begin{equation}
A(t,y)\leq C(1+\norm{t}+\norm{y})^{-N} \label {A}\,,
\end{equation}
which allows to conclude.
Indeed, separating the cases of $\scal{x,y}$ being positive or
negative, we get
\begin{align}
A(t,y)\leq&
\int\frac{e^{-\pi\norm{x-y}^2}}{(1+\norm{t-x}+\norm{y})^N}dx
+\int\frac{e^{-\pi\norm{x+y}^2}}{(1+\norm{t+x}+\norm{y})^N}dx\cr
      =& I_1+I_2.\notag
\end{align}
As $I_2(x,t)=I_1(-x,t)$, it is enough to get a bound for
$I_1$. Now, fix $0<c<1$ and write $B=(1+\norm{t}+\norm{y})$, then
$$
I_1\leq\int_{\norm{x-y}>cB}{e^{-\pi\norm{x-y}^2}dx}
+\int_{\norm{x-y}\leq cB}{\frac{e^{-\pi
\norm{x-y}^2}}{(1+\norm{t-x}
+\norm{y})^N}}dx.
$$
We conclude directly for the first integral.
For the second one, it is sufficient to note that, if $\norm{x-y}\leq c(1+\norm{t}+\norm{y})$, then
\begin{align}
1+\norm{t-x}+\norm{y}\geq&
1+\frac 12\norm{t}+\frac 12\norm{y}-\frac 12\norm{x-y}\cr
\geq& \frac {(1-c)}{2}(1+\norm{t}+\norm{y}).\notag
\end{align}
This completes the proof of (\ref{A}) and (\ref {step2p2}).

Let us finally prove  (\ref {step2p3}).
Fix $c>2$. Then the left hand side is bounded by
\begin{align}
\int_{\norm{x}\leq R}|g(x)|&
\left(\int_{\norm{y}>cR}|\widehat{g}(y)|
e^{2\pi\abs{\scal{x,y}}}dy
+\int_{\norm{y}<cR}|\widehat{g}(y)|e^{2\pi
\abs{\scal{x,y}}}dy\right)dx\cr
&\leq\int_{\norm{x}\leq R}|g(x)|
\left(\int_{\norm{y}>cR}Ce^{-(\pi-2\frac{\pi}{c})\norm{y}^2}dy
+\int_{\norm{y}<cR}|\widehat{g}(y)|e^{2\pi
\abs{\scal{x,y}}}dy\right)dx\cr
&\leq
K\norm{g}_{L^1}+\int_{\norm{x}\leq R}\int_{\norm{y}<cR}|g(x)||\widehat{g}(y)|e^{2\pi
\abs{\scal{x,y}}}dxdy.\notag
\end{align}

Then, if we multiply and divide by
$(1+\norm{x}+\norm{y})^N$ in the integral of right side, we get
the required inequality (\ref {step2p3}).
This completes the proof of the claim.\hfill$\diamond$

\vskip3pt
\noindent{\it Third step.} {\sl The function $g$ admits an holomorphic extension
to $\C^d$ that is of order $2$. Moreover, there exists a polynomial $R$
such that for all $z\in\C^d$, $g(z)g(iz)=R(z)$.}
\vskip3pt
It follows from (\ref{step2p1}) and Fourier inversion that
$g$ admits an holomorphic extension to $\C^d$ which we again denote by $g$.
Moreover,
$$
\abs{g(z)}\leq Ce^{\pi \norm {z}^2}\,,
$$
with $C$ the $L^1$ norm of $\widehat{g}$. It follows that
$g$ is of order $2$.
On the other hand, for all $x\in \R^d$ and $e^{i\theta}$ of modulus $1$,
\begin{equation}
\vert g(e^{i\theta}x)\vert \leq \int_{\R^d}\abs{\widehat{g}}(y)e^{2\pi
\abs{\scal{x,y}}}dy.
\label{complexe}\end{equation}
Let us now define a new function $G$ on $\C^d$ by\,:
$$
G\,:z\rightarrow\int_{0}^{z_1}\dots\int_{0}^{z_d}g(u)g(iu)du.
$$
As $g$ is entire of order $2$, so is $G$. By differentiation of $G$,
the proof of this step is
complete once we show that $G$ is a polynomial.

To do so, we will use (\ref {step2p3}) and an elementary variant
of Phragm\`en-Lindelh\"of's principle which we recall here, and
which may be found in \cite{GS2}\,: {\sl let $\phi$ be an entire
function of order 2 in the complex plane and let $\alpha\in
]0,\pi/2[$; assume that $|\phi(z)|$ is bounded by $C(1+|z|)^N$ on
the boundary of some angular sector
$\{re^{i\beta}\,:r\geq0,\beta_{0}\leq\beta\leq
\beta_{0}+\alpha\}$. Then the same bound is valid inside the
angular sector (when replacing $C$ by $2^NC$).}

Let us fix a  vector $\xi\in\R^d$ and define the function $G_\xi$
on $\C$ by $G_\xi(z)=G(z\xi)$. Then $G_\xi$ is an entire function
of order 2 which has polynomial growth on $\R$ and on $i\R$ by
(\ref{complexe}) and (\ref {step2p3}). We cannot directly  apply
 Phragm\`en-Lindelh\"of's principle since we are not allowed to do so on  angular sectors of
angle $\pi/2$. But, to prove that $G$ has polynomial growth in the
first quadrant, it is sufficient to prove  uniform estimates of
this type inside all angular sectors
$\{re^{i\beta}\,:r\geq0,0<\beta_{0}\leq\beta\leq \pi/2\}$.
Moreover, it is sufficient to have uniform estimates for the
functions $G_\xi^{(\alpha)}(z)=G^{(\alpha)}(z\xi)$, with
$0<\alpha< \beta_0$, and $$
G^{(\alpha)}(z)=\int_{0}^{z_1}\dots\int_{0}^{z_d}g(e^{-i\alpha}u)g(iu)du.
$$ $G_\xi^{(\alpha)}$ clearly has  polynomial growth on $e^{
i\alpha}\R$ and on $i\R$, that is $$ G_\xi^{(\alpha)}(z)\leq
C(1+\abs{x}\norm{\xi})^N. $$ The constant $C$, which comes from
the constant in (\ref {step2p3}),  is independent of $\alpha$. The
same estimate is valid inside the angular sector by the
Phragm\`en-Lindelh\"of's principle, and extends to $G_\xi$, which
we wanted to prove.

Proceeding in an analogous way in the three other quadrants, we prove that
$G_\xi$ is an entire function with polynomial growth of order $N$,
so a polynomial of degree $\leq N$.
Let us now write
$$
G_\xi(z)=a_{0}(\xi)+\dots+a_N(\xi)z^N.
$$
Then
$$
a_j(\xi)=\frac{1}{j!}\frac{d^j}{dz^j}\bigl(G(z\xi)\bigr)\big\vert _{z=0}
$$
shows that $a_j$ is  a homogeneous polynomial of degree $j$ on $\R^d$.

The entire function $G$, which  is a polynomial on $\R^d$, is a
polynomial. Finally,
\begin{equation}
g(z)g(iz)=R(z),\label{prod}
\end{equation}
where $R$ is a polynomial
and the proof of this step is
complete.

\vskip3pt
\noindent{\it Step 4.} {\sl A lemma about entire functions of several variables.}
\vskip3pt

We are now lead to solving the equation (\ref{prod}),
where $g$ is an entire function of order $2$ of $d$ variables and $R$
is a given polynomial. It is certainly well known that such functions $g$ can be
written as $P(z)e^{Q(z)}$, with $Q(z)$ a  polynomial of
degree at most $2$. Moreover, the equation implies that
$Q(z)+Q(iz)=0$, so that $Q$ is homogeneous of degree $2$. So we have
completed the proof, up to the study of the equation (\ref {prod}).
Since we did not find a simple reference for it, we include the proof
of the next lemma, which is a little more general than what we need
above.

\begin{lemma}
Let $\ffi$ be an entire function of order $2$ on $\C^d$ such that, on
every complex line, either $\ffi$ is identically $0$ or it has at
most $N$ zeros. Then, there exists a polynomial $P$ with degree
at most $N$ and a polynomial $Q$ with degree
at most $2$ such that $\ffi(z)=P(z)e^{Q(z)}$.
\label{lempol}
\end{lemma}

\begin{proof}[Proof]  Without loss of generality, we may assume that
$\ffi(0)$ does not vanish.
Then, for $z \in \C^d$, $\ffi_z(t)=\ffi(tz)$ is a non-zero entire
function of order $2$ that has at most
$N$ zeros. By Hadamard's factorisation theorem,
for every $z\in \C^d$,
there exists a polynomial $P_{z}$ with $\deg(P_z)\leq
N$ and $\alpha(z),\beta(z)\in\C$ such that\,:
$$
\ffi(tz)=P_{z}(t)e^{\alpha(z)t+\beta(z)t^2},
$$
with
$$
P_{z}(t)=a_0(z)+\dots+a_N(z)t^N\,.
$$
From the uniqueness in  Hadamard's theorem, we easily see that the
functions $\alpha$ and $\beta$ are homogeneous of degree $1$ and $2$
respectively, and $a_{j}$ is of degree $j$. We may assume that
 $a_N(z)$ is non identically zero.
We have
\begin{equation}
\ffi(tz)\ffi(-tz)\ffi(itz)\ffi(-itz)=P_{z}(t)P_{z}(-t)P_{z}(it)P_{z}(-it).
\label{lempolP1}
\end{equation}
Differentiating (\ref{lempolP1}) $4N$ times with respect to $t$,
and then taking $t=0$, we get that $a_N(z)^4$ is a homogeneous
polynomial of degree $4N$. Now, $$
\frac{d}{dt}\bigl(\ffi(tz)\bigr)\ffi(-tz)\ffi(itz)\ffi(-itz)
=\left(\frac{d}{dt}P_z(t)+
(2\beta(z)t+\alpha(z))P_z(t)\right)P_z(it)P_z(-it), $$ and
differentiating $4N+1$ times with respect to $t$ at $t=0$, we get
that $a_N(z)^4\beta(z)$ is holomorphic. Thus $\beta(z)$ is also
holomorphic, and so  a homogeneous polynomial of degree $2$. An
analogous proof allows to conclude that $\alpha(z)$ is a
homogeneous polynomial of degree $1$. Define
$Q(z)=\beta(z)+\alpha(z)$ and $P(z)=\ffi(z)e^{-Q(z)}$. We know
that $P$ is holomorphic, and we have to prove that it is a
polynomial. Then $$ P(tz)=P_z(t)=a_0(z)+\dots+a_N(z)t^N. $$ In
particular, $$
a_j(z)=\frac{1}{j!}\left.\frac{d^j}{dt^j}\bigl(P(tz)\bigr)\right|_{t=0}
$$ is a holomorphic function, thus a homogeneous polynomial of
degree $j$. It follows that $P(z)=a_0(z)+\dots+a_N(z)$ is a
polynomial of degree $\leq N$, which we wanted to prove to
conclude for this step. \epf We have also completed the proof of
Proposition \ref{BHT2}. Indeed, $g$ has the required form thus so
has $f$. \epf

Let us make a few comments on the proof. Step 3 is very much
inspired, with simplifications, from the proof of H\"ormander
\cite {HB}. Step 2 is not contained in \cite {HB}, and simplifies
greatly the proof, even for $N=0$. Further, it is easy to see that
if the function $\ffi$ in Lemma \ref{lempol} is of order less than
$k$, then $\ffi(z)=P(z)e^{Q(z)}$, with $P$ a polynomial  of degree
at most $N$ and $Q$ a polynomial of degree at most $k$.

As we said in the introduction, the reader will find separately the
proof of Beurling-H\"ormander's Theorem ($N=0$, $d=1$) in the Appendix.

\section{Applications to other uncertainty principles.}
\label{corollaires}

Let us first mention the  following immediate corollary of Theorem
\ref{BHT}\,:

\begin{corollary}
\label{HBP1}
Let $f\in L^2(\R^d)$.
\begin{itemize}
\item[$(i)$]
If
\begin{equation}\dst
\iint_{\R^d\times\R^d}
\frac{|f(x)||\widehat{f}(y)|}{(1+\norm{x}+\norm{y})^N}
e^{2\pi\sum_{i=1}^d\abs{x_iy_i}}dxdy<+\infty
\label{propcor1}
\end{equation}
then $f(x)=P(x)e^{-\sum_{i=1}^d \beta_i x_i^2}$ with $P$  a
polynomial and $\beta_i>0$ for $i=1,\ldots,d$.
\item[$(ii)$] If
\begin{equation}
\iint_{\R^d\times\R^d}
\frac{|f(x)||\widehat{f}(y)|}{(1+\norm{x}+\norm{y})^N}
e^{2\pi\norm{x}\norm{y}}dxdy<+\infty
\label{propcor2}
\end{equation}
then $f(x)=P(x)e^{-\beta \norm{x}^2}$ with $P$ a polynomial and
$\beta>0$.
\end{itemize}
\end{corollary}

\begin{proof}[Proof] It is enough to see that Conditions
(\ref{propcor1}) and (\ref{propcor2}) are stronger than  Condition
(\ref{condBHT}) of Theorem \ref{BHT}. Thus
$f(x)=P(x)e^{-\scal{Ax,x}}$ for some positive definite matrix $A$.
A direct computation then shows that the form of the matrix $A$
imposed by Conditions (\ref{propcor1}) and (\ref{propcor2}) are
respectively $A$ diagonal and $A=\beta I$.
\end{proof}
The next proposition, which implies the   Cowling--Price theorem
in one dimension,  follows at once from the last case.

\begin{proposition}
Let $N\geq 0$. Assume that $f\in L^2(\R^d)$ satisfies $$
\int_{\R^d}\frac{|f(x)|e^{\pi
a\norm{x}^2}}{(1+\norm{x})^N}dx<+\infty \hspace{2ex}\mathrm{and}
\hspace{2ex} \int_{\R^d}\frac{|\widehat{f}(y)|e^{\pi
b\norm{x}^2}}{(1+\norm{y})^N}dy<+\infty \;.$$ Then,  $ab>1$
implies $f=0$. If $ab=1$, then $f(x)=P(x)e^{-\pi a\norm{x}^2}$ for
some polynomial $P$. \label{Cowl}\end{proposition}

Another remark is that $L^1$ norms may be replaced by $L^p$ norms
for $1\leq p\leq \infty$ in the statement of Theorem  \ref{BHT}. In
particular, we get  Corollary \ref{HBP3}  as well as modifications
of the usual Hardy's theorem that we
state now.

\begin{proposition}[Hardy type]
Let $f\in L^2(\R^d)$ be such that $$ |f(x)|\leq
C(1+\norm{x})^Ne^{-\pi\scal{Ax,x}}\hspace{2ex} \mathrm{and}
\hspace{2ex}|\widehat{f}(y)|\leq
C(1+\norm{y})^Ne^{-\pi\scal{By,y}}\,, $$ where $A$ and $B$ are two
real symmetric positive definite matrices. We have the following
implications.
\begin{enumerate}
\item \label{hf1}
If the matrix $B-A^{-1}$ is semi-positive definite and non zero,
then $f=0$.
\item \label{hf2}
If $B=A^{-1}$, then $f(x)=P(x)e^{-\pi\scal{Ax,x}}$ where $P$ is a
polynomial of degree $\leq N$.
\item \label{hf3}
Else, there is a dense subspace of functions satisfying these estimates.
\end{enumerate}
\label{hardyfour}
\end{proposition}

\bpf It is enough to notice that $$
\scal{Ax,x}+\scal{A^{-1}y,y}\geq 2\scal{x,y} $$ and to apply
Theorem \ref{BHT} to obtain the second point. The two other cases
are easy consequences of case \ref{hf2}.
\end{proof}
%Another form of Hardy's theorem states that the integrals
%of $f$ and $\widehat{f}$ outside balls of radius $R$ cannot both
%decrease too fast with respect to $R\to+\infty$\,:
%
%\begin{proposition}
%\label{Hardy integral}
%Let $N\geq0$ be an integer and let $f\in L^{2}(\R^d)$. Assume that
%there is a constant $C$ such that for all $R>0$,
%$$\int_{\norm{x}>R}|f(x)|dx\leq C(1+R)^Ne^{-aR^2}
%\hspace{3ex}\mathrm{and} \hspace{3ex}
% \int_{\norm{y}>R}|\widehat{f}(y)|dy\leq
%C(1+R)^Ne^{-bR^2}
%$$
%with $ab=\pi^2$, then
%$f(x)=P(x)e^{-a\norm{x}^2}$ for some polynomial $P$ of degree at most $N+1$.
%%$N-d+2$.
%\end{proposition}
%
%\bpf To prove this, it is sufficient
%to show that the two integrals
%$$
%\int\frac{|f(x)|}{(1+\norm{x}^2)^m}e^{a \norm{x}^2}dx \hspace{2ex}
%\mathrm{and} \hspace{2ex}
%\int\frac{|\widehat{f}(y)|}{(1+\norm{y}^2)^m}e^{b \norm{y}^2}dy
%$$
%are finite for some integer $m$, which is done by the computation of
%$d$ integrations by parts, and then to apply Proposition \ref{Cowl}.
%It is elementary
%to see that the condition on the degree is necessary and sufficient.
%\epf

\begin{proof}[Proof of Theorem \ref{Cowl1}]
For $d=1$, this is exactly Proposition \ref{Cowl}. For $d>1$, we
proceed as in \cite{BaSw} to reduce to the one-dimensional case.
For almost every $x'=( x_{2}, \cdots x_{d})$, the function
$f_{x'}$ defined by $f_{x'}(x_{1})=f(x_{1},x')$ is in $L^2(\R)$
and has as Fourier transform the function $$ y_{1}\mapsto
\int_{\R^{d-1}}\widehat f(y_{1}, y')e^{-2\pi i\scal{x',y'}}dy'\;.
$$ So, for almost every $x'$, the function $f_{x'}$ satisfies the
assumptions of Proposition \ref{Cowl}, and $e^{\pi
a\abs{x_{1}}^2}\!\! f(x_{1},x')$ is a polynomial of degree at most
$N-2$ in the $x_{1}$ variable. The same is valid in each variable,
which allows to conclude.
\end{proof}
Let us remark that, as in \cite{BaSw} for the case $N=0$, it is
possible to weaken the assumption when the conclusion is that $f$
vanishes. We have the immediate corollary\,:

\begin{corollary}
Assume that $f\in L^2(\R^d)$ satisfies $$
\int_{\R^d}|f(x)|\frac{e^{\pi
a\abs{x_{j}}^2}}{(1+\abs{x_{j}})^M}dx<+\infty
\hspace{2ex}\mathrm{and} \hspace{2ex}
\int_{\R^d}|\widehat{f}(y)|\frac{e^{\pi
b\abs{y_{j}}^2}}{(1+\abs{y_{j}})^N} dy<+\infty $$ for some
$j=1,\cdots,d$ and for some positive constants $a$ and $b$ with
$ab=1$. If $\min\{M,N\}=1$, then $f$ is identically
$0$.\label{Cowl3}
\end{corollary}

\begin{proof}[Proof of Theorem \ref{GelShi}]
It is sufficient to consider the one-dimensional case: for the
first statement we conclude  the general case from the
one-dimensional one as before, and to find a dense subset of
functions we use tensorization. So, let us first assume  that
$f\in L^2(\R)$ is such that $$
\int_{\R}|f(x)|e^{2\pi\frac{a^p}{p}\abs{x}^p}dx<+\infty
\hspace{2ex}\mathrm{and} \hspace{2ex}
\int_{\R}|\widehat{f}(y)|e^{2\pi\frac{b^q}{q}\abs{y}^q}dy<+\infty
$$ for some positive constants $a$ and $b$, with $f$ not
identically zero. Here $1<p\leq 2$, and $q$ is the conjugate
exponent. It follows from the second inequality that $f$ extends
to an entire function, which satisfies the inequality $$
\abs{f(x+iy)}\leq C e^{\frac{2\pi}{p}\abs{\frac{y}{b}}^p}\;. $$
Moreover, the same inequality is valid when $f$ is replaced by its
even part, or its odd part. Such a function may be written as
$g(z^2)$, or $zg(z^2)$, with $g$ an entire function. One of them
is non zero, and satisfies $$ \abs{g(x+iy)}\leq C
e^{2\pi\frac{b^{-p}}{p}\abs{z}^{\frac{p}{2}}}
\hspace{2ex}\mathrm{and} \hspace{2ex}
\int_{\R}|g(x)|e^{2\pi\frac{a^p}{p}\abs{x}^\frac{p}{2}}dx<+\infty\;.
$$ In the second inequality, $a$ has effectively been replaced by
an arbitrarily close smaller constant, which we write $a$ again
for simplification. We then consider the function
 $$G(z)=\int_{0}^z g(u)du-\int _{0}^\infty g(u)du\;. $$
 Then $G$ is a non zero entire function of order $ \frac p2$, and, moreover,
 $$\abs{G(x+iy)}\leq C\abs{z} e^{2\pi\frac{b^{-p}}{p}\abs{z}^{\frac
 p2}}
 \hspace{2ex}
\mbox{and, for positive $x$,} \hspace{2ex} |G(x)|\leq
Ce^{-2\pi\frac{a^p}{p}\abs{x}^\frac p2}\;.$$ Let us consider the
indicator function of $G$, given by
$$h_{G}(\theta)=\limsup_{r\rightarrow \infty}\frac {\ln \abs{G(r
e^{i\theta})}}{r^\frac p2}$$ (see \cite {L} for this topic). Then
there exists a positive measure on the unit circle, $d\mu$, such
that the $2\pi$-periodic function $h_{G}$ is the convolution of
$d\mu$ with the continuous function $\phi_{p}$ defined by
$\phi_{p}(\theta)=\cos(\frac{p\theta}{2})$ on $(0, 2\pi)$. An
elementary computation shows that
$$\abs{\cos(\frac{p\pi}{2})}\phi_{p}(\theta)
+\phi_{p}(\theta+\pi)\geq 0\;.$$ It implies that
$$\abs{\cos(\frac{p\pi}{2})}h_{G}(\theta) +h_{G}(\theta+\pi)\geq
0\;,$$ and, in particular,
$\abs{\cos(\frac{p\pi}{2})}h_{G}(\pi)+h_{G}(0)\geq 0$. But
$h_{G}(\pi)\leq 2\pi\frac{b^{-p}}{p}$, while $h_{G}(0)\leq
-2\pi\frac{a^p}{p}$. The inequality $ab\leq
\abs{\cos(\frac{p\pi}{2})}^\frac 1p$ follows at once.

Let us now prove that, for $ab<\abs{\cos(\frac{p\pi}{2})}^\frac
1p$, there is a dense subset of functions $f$ such that $$
|f(x)|\leq Ce^{-2\pi\frac{a^p}{p}\abs{x}^p} \hspace{2ex}
\mathrm{and} \hspace{2ex} |\widehat{f}(y)|\leq C
e^{-2\pi\frac{b^q}{q}\abs{y}^q}\;. $$ Since this set of functions
is stable under multiplications by $e^{2\pi iy_{0}x}$ (changing
$b$ into an arbitrarily close smaller constant), we see
immediately  that it is dense, unless it reduces to $0$. Indeed,
if $f$ is such a non zero function and $\phi\in L^2(\R)$ is
orthogonal to all functions $e^{2\pi iy_{0}x}f$, then $\overline
f\phi$ is identically $0$. Since $f$ is analytic, it means that
$\phi$ is $0$.

Moreover,  using a dilation if necessary, we can restrict to the
case when $a<\abs{\cos(\frac{p\pi}{2})}^{\frac 1p}$ and $b<1$. We
shall prove a slightly stronger result, the fact that $E$ is not
reduced to $0$, where $E$ is the class of functions for which, for
every $a<1$, there exists a constant $C_{a}$ such that
\begin{equation}
|f(x)|\leq C_{a}e^{-2\pi\frac{a^p}{p}\abs{\cos(\frac{p\pi}{2})}\abs{x}^p}
\hspace{2ex}\mathrm{and} \hspace{2ex}
|\widehat{f}(y)|\leq C_{a}e^{-2\pi\frac{a^q}{q}\abs{y}^q}\;.
\label{four-exp}
\end{equation}
We claim that this follows easily  from the following lemma:

\begin{lemma}
The function $f$ is in $E$ if and only if it extends to an entire
function of order $p$, and, for every $a<1$, there exists a
constant $C_{a}$ such that
\begin{equation}
|f(x)|
\leq C_{a}e^{-2\pi\frac{a^p}{p}\abs{\cos(\frac{p\pi}{2})}\abs{x}^p}
\hspace{2ex}\mathrm{and} \hspace{2ex}
|{f}(iy)|\leq C_{a}e^{2\pi\frac{a^{-p}}{p}\abs{y}^p}\;.
\label{compl-plane}\end{equation}
\end{lemma}

Indeed, let us take the lemma for granted and finish the proof.
We know that there exists non zero entire functions of order $p$ which satisfy
these two conditions (see  \cite{GS2} or \cite{L}). So, it
follows from the lemma that there exists non zero functions in $E$.
\begin{proof}[Proof of the lemma]
One implication (not the one that we need) is trivial. So let us
assume that $f$ is an entire function which satisfies
(\ref{compl-plane}).  We adapt the proof from \cite{GS2} with the
required modifications to have a precise control of the constants.
We first give precise estimates of the function inside each of the
quadrants of the complex plane which are delimited by the
coordinate axes. It follows from Phragm\`en-Lindel\"of principle
and the assumption that, for every $\varepsilon>0$ there exists a
constant $C_{ \varepsilon}$ such that $\theta\in[0, \pi/2]$ $$
\abs{ f(re^{i\theta})}\leq C_{\varepsilon}
e^{\frac{2\pi}{p}\cos(p(\theta -\frac{\pi}{2}))r^p}e^{\varepsilon
r^p}\;.$$ From the convexity of the function $\theta\mapsto -\log
(\cos \theta)$, we see  that $\cos(\frac{p(\theta -\pi}{2}))\leq
(\sin \theta)^p$. Combining the corresponding inequality for $f$
with a simple estimate of $f$ in an angular sector around the
$x$-axis, and doing the same in the three other quadrants, we
deduce the existence of constants $\alpha>1$ and $\beta>0$ such
that,   for every $a>1$, there exists a constant $C_{a}$, such
that,
\begin{equation}
\abs{ f(z)}\leq C_{a}e^{a\frac{2\pi}{p}\abs{y}^p}
\hspace{3ex}\mathrm{and,\ for}\ \abs{x}\geq \alpha \abs{y}\,,
\hspace{2ex}
\abs{ f(z)}\leq C_{a} e^{-\beta\abs{x}^p}\,.
\label{fin}
\end{equation}
Let us now deduce from these last inequalities the required
exponential decrease of the Fourier transform of $f$. It is
sufficient to prove that, for every $a>1$, there exists a contant
$C_{a}$ such that, for all positive integers $n$ and for all
$y\in\R$, $$
\left(\frac{2\pi}{q}\right)^n\abs{y}^{nq}\abs{\widehat f(y)}\leq
C_{a}a^n\Gamma(n+1)\;. $$ This is equivalent to proving (possibly
by changing the constant $C_{a}$), that, for all positive integers
$n$ and for all $y\in\R$, $$ \abs{y}^{n}\abs{\widehat f(y)}\leq
C_{a}a^n\left(\frac{q}{2\pi}\right)^{\frac nq}\Gamma(\frac nq
+1)\;. $$ This in turn is implied by the inequality $$
\int_{\R}\abs{ f^{(n)}(x)}dx\leq C_{a}a^n
q^{\frac{n}{q}}\Gamma\left(\frac{n}{q}+1\right)
(2\pi)^{\frac{n}{p}}\;, $$ which we now prove using (\ref{fin}).
To estimate $\abs{ f^{(n)}(x)}$, we use the Cauchy inequality
related to the circle centered at $x$ and with radius
$R=\left(\frac{n}{2\pi}\right) ^{\frac 1p}$. For $\abs{x}\leq
2\alpha R$, we use the first part of (\ref{fin}), to get the
inequality $$ \abs{ f^{(n)}(x)}\leq C_{a}\frac {\Gamma
(n+1)}{R^{n}} e^{a\frac{2\pi}{p}R^p}\leq C'_{a}a^n (2\pi)^\frac np
e^\frac np \frac {\Gamma (n+1)}{n^{\frac np}}\;. $$ Then
Stirling's formula for the asymptotics of $\Gamma$ allows to
conclude for the integral on the interval $\abs{x}\leq 2\alpha R$.
Outside, we get the estimate $$ \abs{ f^{(n)}(x)}\leq
C_{a}(2\pi)^\frac np \frac {\Gamma (n+1)}{n^{\frac
np}}e^{-\beta'\abs{x}^p}\;. $$ We conclude as before using the
Stirling's formula. \epf This finishes the proof of Theorem
\ref{GelShi}. \epf It leaves completely open the study of the
critical case, $ab= \abs{\cos(\frac{p\pi}{2})}^\frac 1p$, except
for the case when $p=2$, for which the previous theorems give
precise informations.

\section{Properties of the ambiguity function} \label{propamb}

For sake of self-containedness, let us recall here a few
properties of the ambiguity function that we may use in the
sequel. They can be found in \cite{AT}, \cite{Wi}.
\begin{lem}\label{lem0}
Let $u,v$  in $L^2(\R^d)$.
For $a,\omega\in\R^d$, $\lambda>0$ define
$$
S(a)u(t)=u(t-a)\quad,\quad M(\omega)u(t)=e^{2i\pi\omega t}u(t)
\quad,\quad D_\lambda u(t)=\lambda^\frac d2 u(\lambda t)
$$
and recall that we defined $Zu(t)=u(-t)$. Then
\begin{tabbing}
\quad\=$(iii)$ \=$A\bigl(D_\lambda u,D_\lambda v\bigr)(x,y)$\kill
\>$(i)$
\>$A\bigl(S(a)u,S(b)v\bigr)(x,y)=e^{i\pi\scal{a+b,y}}A(u,v)(x+b-a,y)$,\\
\>$(ii)$ \>$A\bigl(M(\omega_1)u,M(\omega_2)v\bigr)(x,y)=
e^{i\pi\scal{\omega_1+\omega_2,x}}A(u,v)(x,y+\omega_1-\omega_2)$,\\
\>$(iii)$ \>$A\bigl(D_\lambda u,D_\lambda v\bigr)(x,y)=A(u,v)\left(\lambda
x,\frac{y}{\lambda}\right)$,\\
\>$(iv)$ \>$A(Zu,Zv)(x,y)=\overline{A(u,v)(x,y)}$,\\
\>$(v)$ \>$A(\widehat u,\widehat v)(x,y)=A(u,v)(y,-x)$ and
$A(u,v)(x,y)=\overline{A(v,u)(-x,-y)}$.\\
\end{tabbing}
\end{lem}
\begin{lemma}
\label{lem1}
Let $u,v\in L^2(\R^d)$. Then $A(u,v)$ is continuous on $\R^{2d}$ and
$A(u,v)\in L^2(\R^{2d})$. Further,
$$
\norm{A(u,v)}_{L^2(\R^{2d})}=
\norm{u}_{L^2(\R^{d})}\norm{v}_{L^2(\R^{d})}
$$
\end{lemma}

\begin{proof}[Proof] This fact is also well known, however to help
the reader  get familiar with our notation, let us recall the
proof of the last assertion (see \cite{AT}, \cite{Wi}). If $u$ and
$v$ are fixed, we will write
\begin{equation}
    h_x(t)=u(t+\frac{x}{2})\overline{v(t-\frac{x}{2})}
\label{hx}\end{equation}
The change of variables
\begin{equation}
  \xi=t-\frac{x}{2} \hspace{4ex}
\mathrm{and} \hspace{4ex} \eta=t+\frac{x}{2}
\label{change}\end{equation}
gives
$$
\iint{\abs{h_x(t)}^2dtdx}
=\iint{\abs{u(\eta)}^2\abs{v(\xi)}^2d\eta d\xi}=\norm{u}_{L^2(\R^{d})}^2
\norm{v}_{L^2(\R^{d})}^2.
$$
In particular, for almost every $x$, the integral with respect to
$t$ is finite, {\it i.e} $h_x\in L^2$. Noticing that
$A(u)(x,y)=\widehat{h_x}(y)$, and using Parseval's formula
we obtain
$$
\int\left(\int\abs{A(u,v)(x,y)}^2dx\right)dy
=\int\left(\int\abs{\widehat{h_x}(y)}^2dy\right)dx
=\norm{u}_{L^2(\R^{d})}^2\norm{v}_{L^2(\R^{d})}^2,
$$
which completes the proof.
\end{proof}
Finally, we will also need the following lemma from \cite{Janote},
\cite{Janssen}\,:

\begin{lemma} Let $u,v,w\in L^2(\R^d)$. Then, for every $x,y\in \R^d$,
$$
\int_{\R^{2d}}A(u,v)(s,t)\overline{A(v,w)(s,t)}
e^{2i\pi(\scal{s,x}+\scal{t,y})}dsdt
=A(u,v)(-y,x)\overline{A(v,w)(-y,x)}.
$$
\label{fouramb}
\end{lemma}

\section{Heisenberg inequality for the ambiguity function.}
\label{heisamb}

We show here that the ambiguity function is subject to sharp
inequalities of Heisenberg type. We give first a directional
version of Heisenberg's inequality (\ref{eqth2}) in the context of
ambiguity functions, with an elementary proof.

\begin{theorem}
\label{th3}
For $u,v\in L^2(\R^d)$, for every
$i=1,\ldots,d$ and every $a,b\in\R$, one has the following inequality\,:
\begin{equation}
\dst\iint_{\R^{2d}}{\abs{x_i-a}^2\abs{A(u,v)(x,y)}^2dxdy}
\dst\iint_{\R^{2d}}{\abs{y_i-b}^2\abs{A(u,v)(x,y)}^2dxdy}
\geq
\frac{\norm{u}_{L^2(\R^d)}^4\norm{v}_{L^2(\R^d)}^4}{4\pi^2}.
\label{ineqth3}
\end{equation}
Moreover equality holds in (\ref{ineqth3}), with $u$ and $v$
non identically $0$, if and only if there exists
$\mu,\nu\in L^2(\R^{d-1})$, $\alpha>0$ and $\beta, \gamma \in\R$
such that
\begin{align}
u(t)=&\mu(t_1,\ldots,t_{i-1},t_{i+1},\ldots,t_d)
e^{2i\pi \beta t_i}e^{-\alpha /2\abs{t_i-\gamma}^2},%\ \mathrm{and}
\notag\\
v(t)=&\nu(t_1,\ldots,t_{i-1},t_{i+1},\ldots,t_d)
e^{2i\pi (b+\beta)t_i}e^{-\alpha /2\abs{t_i-a-\gamma}^2}.\notag
\end{align}
\end{theorem}

To prove the theorem, we will need the following lemma
that has its own interest\,:
\begin{lemma}
\label{lem2}
Let $u,v\in L^2(\R^d)$ be both non identically zero and let
$i=1,\ldots,d$. The following are equivalent\,:
\begin{itemize}
\item[$(a)$]\label{implem21}
$\dst\iint{|x_i|^2\abs{A(u,v)(x,y)}^2dxdy}<+\infty.$

\item[$(b)$]\label{implem22} For all $a\in\R$,
$\dst\iint{|x_i-a|^2\abs{A(u,v)(x,y)}^2dxdy}<+\infty.$

\item[$(c)$]\label{implem23}
$\dst\int{|t_i|^2\abs{u(t)}^2dt}<+\infty$ and
$\dst\int{|t_i|^2\abs{v(t)}^2dt}<+\infty.$
\item[$(d)$]\label{implem24} For all $a,b\in\R$,
$\dst\int{|t_i-a|^2\abs{u(t)}^2dt}<+\infty$ and
$\dst\int{|t_i-b|^2\abs{v(t)}^2dt}<+\infty.$
\end{itemize}
\end{lemma}

\begin{remark} Note that, if $u\in L^2(\R^d)$ and if
$\displaystyle\int{|t_i|^2\abs{u(t)}^2dt}<+\infty$ then
$t_i\abs{u(t)}^2\in L^1(\R^d)$.
\end{remark}

\begin{proof}[Proof]
 Let us first remark that $(a)\Leftrightarrow
(b)$ and $(c)\Leftrightarrow (d)$. Indeed, it is sufficient to use
the triangle inequality $$ |x|^2\leq 2(|x-a|^2+|a|^2)\,. $$
Moreover, $(b)$ and $(d)$ may be replaced by $(b')$ and $(d')$,
where ``{\it for all}'' has been changed into ``{\it for some}''.
Let us prove that $(a)$ implies $(d')$. With the notations of
Lemma \ref{lem1}, Parseval identity gives\,:
\begin{align}
\iint{|x_i|^2\abs{A(u,v)(x,y)}^2dxdy}
&=\int |x_i|^2\int{\abs{{\widehat{h_x}(y)}}^2dy}dx \cr
&=\int |x_i|^2\int{\abs{h_x(t)}^2dt}dx\cr
&=\iint{|x_i|^2\abs{u(t+\frac{x}{2})\overline{v(t-\frac{x}{2})}}^2dtdx}\cr
&=\int\left
(\int{|\eta_i-\xi_i|^2\abs{v(\xi)}^2d\xi}\right
)\abs{u(\eta)}^2 d\eta.
\label{parseval1}
\end{align}
So, if $(a)$ holds, for almost every $\eta$,
$\abs{u(\eta)}^2\int{|\xi_i-\zeta_i|^2\abs{v(\xi)}^2d\xi}
<+\infty$. As we assumed that $u\neq 0$, there exists  $\eta$ such
that $u(\eta)\neq 0$, and the first inequality in $(d')$ holds
with $a=\eta$. Since $u$ and $v$ play the same role, we conclude
for the second part similarly.

Conversely, if $(c)$ holds, the right hand side of (\ref{parseval1})
is finite, and $(a)$ holds also.
\end{proof}

\begin{proof}[Proof of theorem \ref{th3}]
Let us start by proving the inequality.
Set $A=A(u,v)$. We may assume that
$$
\left [\iint{|x_i|^2\abs{A(x,y)}^2dxdy}\right ]
\left [\iint{|y_i|^2\abs{A(x,y)}^2dxdy}\right ]<+\infty
$$
so that both factors are finite, and, by homogeneity,
that $\norm{u}_{L^2(\R^d)}=\norm{v}_{L^2(\R^d)}=1$.
Moreover, replacing $u$ and $v$ by translates of these functions (with
the same translation),
we may assume that
\begin{equation}
    -\int t_i\abs{u(t)}^2dt=\int t_i\abs{v(t)}^2dt=\frac{\overline a}{2}\,.
\label{transl}\end{equation}
Heisenberg's inequality (\ref{eqth2}) applied to
$h_x(t)=u\left(t+\frac{x}{2}\right)
\overline{v\left(t-\frac{x}{2}\right)}$
implies that, for any $b\in\R^d$\,:
   \begin{equation}
\frac{1}{4\pi}\int{\abs{h_x(t)}^2dt}
\leq \left(\int{|t_i|^2\abs{h_x(t)}^2dt}\right)^{1/2}
\left(\int{|y_i-b|^2\abs{A(u,v)(x,y)}^2dy}\right)^{1/2}.
\label{incertitude}
\end{equation}
Integrating this inequality with respect to the $x$-variable and
appealing to Cauchy-Schwarz' inequality, we get\,:
\begin{equation}
\frac{1}{4\pi}=\frac{\norm{u}^2\norm{v}^2}{4\pi}\leq
\left(\iint{|t_i|^2\abs{h_x(t)}^2
dtdx}\right)^{\frac{1}{2}}\left(\iint{|y_i-b|^2\abs{A(u,v)(x,y)}^2
dxdy}\right)^{\frac{1}{2}}.
\label{CS}
\end{equation}
Let us now transform the first factor on the right hand side of
this expression. We write
$$
\abs{t_i}^2=
\left(t_i+\frac{x_i-a}{2}\right)\left(t_i-\frac{x_i-a}{2}\right)
+\frac{|x_i-a|^2}{4}.
$$
The second term which appears is, by Parseval identity, equal to
$$
\frac{1}{4}\iint{|x_i-a|^2\abs{A(x,y)}^2dxdy}\,;
$$
The first term is equal to
\begin{equation}
\iint\left(\eta_i-\frac{a}{2}\right)\left(\xi_i+\frac{a}{2}\right)
\abs{u(\eta)}^2\abs{v(\xi)}^2d\eta d\xi
=-\frac{1}{4}\abs{\overline a-a}^2
%\norm{u}_{L^2(\R^d)}^2\norm{v}_{L^2(\R^d)}^2
\leq0
\label{zero}
\end{equation}
using (\ref{transl}).
Finally, including these results in (\ref{CS}), we obtain the desired
inequality.

\medskip

Assume now that we have equality. Let us remark that, using
properties $(a)$ and $(b)$ of Lemma \ref {lem0}, we may as well assume
that the constants $a$ and $b$ are $0$.
Moreover, up to a same translations in space and frequency, we may again assume
that $-\int
t_i\abs{u(t)}^2dt=\int t_i\abs{v(t)}^2dt=\frac{\overline a}{2}$,
and
$-\int t_i\abs{\ff u(t)}^2dt=\int t_i\abs{\ff v(t)}^2dt=\frac{\overline
b}{2}$.
Let $h_x(t)=u\left(t+\frac{x}{2}\right)
\overline{v\left(t-\frac{x}{2}\right)}$ as before.
Then, to have equality
in (\ref{ineqth3}), we have equality in (\ref{zero}),
{\it i.e.} we have $\overline a=a=0$. Similarly, exchanging the roles
of  the $x$ and $y$ variables,
 we also have $\overline b=b=0$.

We then  have equality in Cauchy-Schwarz's inequality (\ref{CS}).
This implies that $$ x\mapsto\dst\int\abs{t_i}^2\abs{h_x(t)}^2dt\
\mathrm{and}\
x\mapsto\dst\int\abs{y_i}^2\abs{\widehat{h_x}(y)}^2dy $$ are
proportional.

Further, for almost
every $x$, we also  have equality
in Heisenberg's inequality (\ref{incertitude}). From now on, we assume
for simplicity that $i=1$. We then get that, for almost every $x$,
$h_x$ is a Gaussian in the $t_1$ variable\,:
$$
h_x(t)=C(x,t')e^{-\pi\alpha(x)\abs{t_1}^2}.
$$
where $t'=(t_2,\ldots,t_d)\in\R^{d-1}$,
$C(x,t')$ is not identically $0$ and $\alpha(x)>0$ for those $x$
for which $C(x,t')\not=0$.

Let us first prove that $\alpha$ does not depend on $x$. With this
expression of $h_x$, we get that $$
\int\abs{t_1}^2\abs{h_x(t)}^2dt=\norm{C(x,.)}^2_{L^2(\R^{d-1})}
\frac{\sqrt{\pi}}{2}\left(\frac{1}{2\pi\alpha(x)}\right)^{\frac{3}{2}}
$$ whereas (with Parseval identity) $$
\int\abs{y_i}^2\abs{\widehat{h_x}(y)}^2dy=\norm{C(x,.)}^2_{L^2(\R^{d-1})}
\frac{1}{4\pi}\sqrt{\frac{\alpha(x)}{2}}. $$ But these two
functions are proportional only if $\alpha(x)$ is a constant, say
$\alpha(x)=\alpha$.

Taking inverse Fourier transforms, we get that
\begin{equation}
A(u,v)(x,y)=\frac{1}{\sqrt{\alpha}}\widehat C(x,-\widehat y')e^{-\frac{\pi}{\alpha}\abs{y_1}^2},
\label{auxygauss}
\end{equation}
where $\widehat C(x,.)$ is the Fourier transform in $\R^{d-1}$ of $C(x,.)$.
In particular, as $A(u,v)$ is continuous, $\widehat C$ is also continuous.

Further, if one has equality in (\ref{ineqth3}) for $u,v$, one has again
this equality if $u,v$ are replaced by their Fourier transform, as
$A(\widehat u,\widehat v)(x,y)=A(u,v)(y,-x)$. So, we have
\begin{equation}
A(u,v)(x,y)=D(y,x')e^{-\frac{\pi}{\beta}\abs{x_1}^2},
\label{auyxgauss}
\end{equation}
for some function $D$ and some $\beta>0$.
Comparing (\ref{auxygauss}) and (\ref{auyxgauss}), we get that
$$
A(u,v)(x,y)=E(x',y')e^{-\frac{\pi}{\beta}\abs{x_1}^2}
e^{-\frac{\pi}{\alpha}\abs{y_1}^2}
$$
for some function $E\in L^2(\R^{d-1}\times\R^{d-1})$.
Taking the Fourier transform in the $y$ variable, we get
$$
u\left(t+\frac{x}{2}\right)\overline{v\left(t-\frac{x}{2}\right)}
=\widetilde E(x',t') e^{-\frac{\pi}{\beta}\abs{x_1}^2}
e^{-\pi\alpha\abs{t_1}^2}.
$$
So, setting again $\eta=t+\frac{x}{2}$, $\xi=t-\frac{x}{2}$, we get
$$
u(\eta)\overline{v(\xi)}=\widetilde
E\left(\eta'-\xi',\frac{\eta'+\xi'}{2}\right)
e^{-\frac{\pi}{\beta}(\eta_1-\xi_1)^2}
e^{-\frac{\pi\alpha}{4}(\eta_1+\xi_1)^2}.
$$
This is only possible if $\beta=\frac{4}{\alpha}$ and
$u(x)=\mu(\widehat x_i)e^{-\frac{\pi\alpha}{2}\abs{x_i}^2}$,
$v(x)=\nu(\widehat x_i) e^{-\frac{\pi\alpha}{2}\abs{x_i}^2}$
with $\mu,\nu\in L^2(\R^{d-1})$.%\hfill$\Box$
\end{proof}

\medskip
We now translate the last theorem in terms of inequalities on
covariance matrices, as it is classical for inequality (\ref{eqth2}) (see
for instance \cite{DCT}).
\begin{proof}[Proof of Theorem \ref{variance}]
First, it follows from Lemma \ref{lem2} that $X$ and $Y$ have
moments of order $2$. Let us  prove that $X$ and $Y$ are not
correlated. Without loss of generality, we may assume that $$ \int
t_i\abs{u(t)}^2dt=\int t_i\abs{v(t)}^2dt= \int t_i\abs{\ff
u(t)}^2dt=\int t_i\abs{\ff v(t)}^2dt=0 $$ for all $i$, so that $X$
and $Y$ are centered. Let us show that $$ \iint
x_{i}y_{j}\abs{A(u,v)}^2dx\, dy=0\;. $$
 Using Plancherel identity
, we are led to consider $$ \iint x_{i}\left
[\partial_{t_{j}}h_{x}(t)\overline {h_{x}(t)}+ h_{x}(t)\overline
{\partial_{t_{j}}h_{x}(t)}\right ]dxdt\;. $$ Writing
$x_{i}=t+\frac{x_{i}}{2}-(t-\frac{x_{i}}{2})$, and $$
\partial_{t_{j}}h_{x}(t)=u
(t+\frac{x}{2})\overline{\partial_{t_{j}}v(t-\frac{x}{2})}+\partial_{t_{j}}u
(t+\frac{x}{2})\overline { v(t-\frac{x}{2})}, $$ we get eight
terms. For four of them we get directly $0$. The sum of the last
four may be written, after changes of variables, as $$
\int_{\R^d}t_{i}\partial_{t_{j}}[\abs{u(t)}^2]dt\times
\int_{\R^d}\abs{v(t)}^2dt - \int_{\R^d}\abs{u(t)}^2dt
\times\int_{\R^d}t_{i}\partial_{t_{j}}[\abs{v(t)}^2]dt\;.$$ After
an integration by parts, the two terms give $\delta_{ij}$, so that
their difference is $0$.

\medskip
Let us now prove the second assertion. Let $C$ be an automorphism
of $\R^d$. For a function $f\in L^2(\R^d)$, we consider $f_{C}$
the function given by $f_{C}(t)=\abs{\det(C)}^{-\frac
12}f(C^{-1}t)$. Then a simple change of variables shows that the
 probability density of $(CX, ^t C^{-1}Y)$ is
$\abs{A(u_{C},v_{C})}^2$. Eventually changing $u$ and $v$ into
$u_{C}$ and $v_{C}$ we may assume that $V(Y)$ is the identity
matrix. Moreover, we may also assume that $V(X)$ is diagonal. Then
the inequality follows from Theorem  \ref{th3}. Equality holds
only if all eigenvalues are equal to $4\pi^2$, which means that
$u$ and $v$ are Gaussian functions.
\end{proof}

\medskip
From Theorem \ref{variance}, we immediately obtain that
\begin{equation}
    \det (V(X))\times \det (V(Y))\geq (4\pi^2)^{-2d}\;.
\label{det}\end{equation}
 Equality holds only for Gaussian
functions. Another (much less elementary) proof of (\ref{det}) can
be obtained using the entropy inequality of Lieb and the theorem
of Shannon (see \cite{FS}, Section 6).

\medskip
The same inequality
holds for traces instead of determinants. We state it independently.
\begin{corollary}
\label{cor4}
Let $u,v\in L^2(\R^d)$ be both non identically zero and $a,b\in \R^d$. Then
\begin{align}
\dst\iint_{\R^{2d}}{\norm{x-a}^2\abs{A(u,v)(x,y)}^2dxdy}&
\times\dst\iint_{\R^{2d}}{\norm{y-b}^2\abs{A(u,v)(x,y)}^2dxdy}\notag\\
&\geq
\frac{d^2\norm{u}_{L^2(\R^d)}^4\norm{v}_{L^2(\R^d)}^4}{4\pi^2},
\label{ineqcor4}
\end{align}
with equality if and only if there exists $\lambda,\nu\in\C^*$,
$\alpha>0$ and $\beta,\gamma\in\R^d$ such that $$ u(t)=\lambda
e^{2i\pi\scal{\beta,t}}e^{-\alpha /2\|t-\gamma\|^2}\ \ {\mathrm
and }\ \ v(t)=\nu e^{2i\pi\scal{b+\beta,t}}e^{-\alpha
/2\|t-a-\gamma\|^2}\,. $$
\end{corollary}

\section{Uncertainty principles for the ambiguity
function.}
\label{BHamb}

We first prove the following uncertainty principle for the
ambiguity function, which also gives a characterization of Hermite
functions:
\begin{theorem}
Let $u,v\in L^2(\R^d)$ be non identically vanishing. If
\begin{equation}
\iint_{\R^{2d}}
\abs{A(u,v)}^2\frac {e^{\pi\norm{x}^2}}{(1+\norm{x})^N}
dxdy<+\infty \hspace{2ex}\mathrm{and} \hspace{2ex}
\iint_{\R^{2d}}
\abs{A(u,v)}^2\frac {e^{\pi\norm{y}^2}}{(1+\norm{y})^N}
dxdy<+\infty
\label{HBA0}
\end{equation}
for all $j=1, \cdots,d$, then there exists $a,w\in\R^d$ such that
both $u$ and $v$ are  of the form $$P(x)e^{2i\pi
\scal{w,x}}e^{\!-\pi\norm{x-a}^2}\,,$$
 with $P$ a polynomial.
\label{thBHamb0}
\end{theorem}
Let us first remark that when $u=v$, since the Fourier transform
of $\abs{A(u,u)}^2$ taken at $(x,y)$ is equal to
$\abs{A(u,u)}^2(y,-x)$, the result follows immediately from
Theorem \ref{Cowl1}. Our aim is to prove that it is also valid in
the general case. We first prove a weaker statement.
\begin{proposition}
Let $u, v\in L^2(\R^d)$ be non identically vanishing. If
\begin{equation}
\iint_{\R^d\times\R^d}
\frac{\abs{A(u,v)(x,y)}^2}{(1+\norm{x}+\norm{y})^N}
e^{\pi(\norm{x}^2+\norm{y}^2)}dxdy<+\infty\,,
\label{HBA}
\end{equation}
then there exists $a,w\in\R^d$ such that both $u$ and $v$ are  of
the form $P(x)e^{2i\pi \scal{w,x}}e^{-\pi\norm{x-a}^2 }$ with $P$
a polynomial. \label{thBHamb}
\end{proposition}

Before starting the proof of Proposition \ref{thBHamb}, let us state two
 lemmas. The first one is elementary and well known.

\begin{lem}
Let $u, v\in L^2(\R^d)$  be non identically vanishing. Then $$
u(x)=P(x)e^{2i\pi \scal{\alpha,x}}e^{-\pi \norm{x-a}^2} \
\mathrm{and }  \ v(x)=Q(x)e^{2i\pi
\scal{\alpha,x}}e^{-\pi\norm{x-a}^2} $$ with $P,Q$ polynomials and
$a,\alpha\in\R^d$ if and only if there is a polynomial $R$ such
that $$ A(u,v)(x,y)=R(x,y)e^{2i\pi(\scal{\alpha,x}+\scal{a,y})}
e^{-\frac{\pi}{2}(\norm{x}^2+\norm{y}^2)}. $$ Moreover,
$\deg(R)=\deg(P)+\deg(Q)$. \label{polyn}\end{lem}

\begin{lem}
Assume that $u, v\in L^2(\R^d)$, with $u(x)=P(x)e^{2i\pi
\scal{a,x}}e^{-\pi \norm{x}^2}$, where $a\in \C^d$. Then the
function $A(u,v)$ can be extended to an entire function on
$\C^{2d}$. \label{entire-amb}
\end{lem}

\bpf For $z,\zeta \in \C^d$, we  note $\scal{z,\zeta}=\sum
z_{i}\zeta_{i}$. Then
\begin{equation}
A(u,v)(x,y)= e^{i\pi\scal{x,(y+2a)}} \times\int P(t+x)\overline {v(t)}e^{-\pi
\norm{t+x}^2}e^{2i\pi \scal{t,(a+y)}}dt\,.
\label{compl}
\end{equation}
This clearly makes sense for $x,y\in \C^d$, and defines an entire
function. \epf

\begin{proof}[Proof of Proposition \ref{thBHamb} ]
By homogeneity, we may assume that
$\norm{u}_{L^2}=\norm{v}_{L^2}=1$. For each $w\in L^2(\R^d)$, we
consider the function $$ g_w=A(u,v)\overline{A(v,w)}\,. $$ By
Lemma \ref{fouramb}, the Fourier transform of $g_{w}$ is given by
$\ff g_w(x,y)=g_w(y,-x)$.

\noindent{\it First step.} {\sl There exists a polynomial $P$ such that}
\begin{equation}
    A(u,v)(x,y)\overline{A(v,u)}(x,y)=P(x,y)e^{-\pi(\norm{x}^2+\norm{y}^2)}\,.
\label{prod-A}
\end{equation}

\bpf We consider here the function $g_u$. As $g_u$ is (up to a
change of variable) its own Fourier transform, by Proposition
\ref{Cowl}, it is sufficient to prove that
\begin{equation}
\int
\abs{g_{u}(x,y)}\frac{e^{\pi(\norm{x}^2+\norm{y}^2)}}{(1+\norm{x}+\norm
{y})^N}dxdy<\infty\,.
\label{finite-A}
\end{equation}
It follows immediately from the assumption on $A(u,v)$, using
Cauchy-Schwarz inequality and the fact that $\overline
{A(v,u)(x,y)}=A(u,v)(-x,-y)$. \epf To complete the proof of the
proposition, it is sufficient to prove that $A(u,v)$ extends to an
entire function of order $2$. Indeed, Lemma \ref{lempol} then
implies that $$ A(u,v)(x,y)=R(x,y)e^{Q(x,y)}\,, $$ where $R$ is a
polynomial and $Q$ a polynomial of degree at most $2$. But, as $$
A(u,v)(x,y)A(u,v)(-x,-y)=A(u,v)(x,y)\overline
{A(v,u)(x,y)}=P(x,y)e^{-\pi(\norm{x}^2+\norm{y}^2)}, $$ we get
$Q(x,y)=\scal{\beta,x}+\scal{\gamma,y}-\frac
\pi2(\norm{x}^2+\norm{y}^2)$ for some constants
$\beta,\gamma\in\C^d$. Next, Condition (\ref{HBA}) implies that
$\beta,\gamma$ are purely imaginary,
$\beta=2i\pi\alpha,\gamma=2i\pi a$ with $a,\alpha\in\R^d$ so that
$$ A(u,v)(x,y)=R(x,y)e^{2i\pi(\scal{\alpha,x}+\scal{a,y})}
e^{-\frac{\pi}{2}(\norm{x}^2+\norm{y}^2)}\,, $$ with $R$ a
polynomial.
 Lemma \ref{polyn} allows to conclude. So, we have
finished the proof once we have proved the second step.
\medskip

\noindent{\it Second step.} {\sl The function $A(u,v)$ extends to
an entire function of order $2$.}

\bpf To prove this, we first show that, for each $w\in L^2(\R^d)$,
the function $g_{w}$ extends to an entire function of order $2$.
Since, up to a change of variable, $g_{w}$ coincides with its
Fourier transform, it is sufficient to show that $$ \int
|g_w(x,y)|e^{\frac{\pi}{4}(\norm{x}^2+\norm{y}^2)}dxdy<\infty\,.
$$ This last inequality follows from the fact that
$\abs{g_w}\leq\abs{A(u,v)}\norm{w}_{L^2}$, and from the assumption
on $A(u,v)$. We get the estimate
\begin{equation}
    |g_{w}(z,\zeta)|\leq C \norm{w} e^{4\pi
    (\norm{z}^2+\norm{\zeta}^2)}\,,
\label{maj-g}\end{equation} for all $(z,\zeta)\in \C^d\times
\C^d$, with a constant $C$ which does not depend on $w$. In order
to conclude  this step,  it is sufficient  to show  that there
exists a constant $C$ such that, for each $(z, \zeta) \in
\C^d\times \C^d$, we may find $w_{z,\zeta}$ which is of the form
required  in Lemma \ref{entire-amb},
  such that
\begin{equation}
    \abs{A(w_{z,\zeta}, v)(z,\zeta)}\geq
    C^{-1}e^{-C(\norm{z}^2+\norm{\zeta}^2)}\,,
\label {below}\end{equation}
and
\begin{equation}
    \norm{w_{z,\zeta}}\leq C^{-1}e^{C(\norm{z}^2+\norm{\zeta}^2)}\,.
\label {above}\end{equation} By density of the Hermite functions
we can choose a polynomial $P_0$ such that $$ \int
P_0(t)\overline{v(t)}e^{-\pi\norm{t}^2}dt=1\,. $$ We then define
$w_{z,\zeta}$ by $$ w_{z, \zeta}(t)=P_0(t-z) e^{2\pi
\scal{z-i\zeta,t}}e^{-\pi \norm{t}^2}\;. $$ It follows from the
choice of $P_{0}$ that $$ A(w_{z,\zeta},v)(z,\zeta)=e^{\pi
(\scal{z, z}-i\scal{z,\zeta})} \,. $$ Then (\ref{below}) and
(\ref{above}) follow from direct computations. Finally, since
$A(u,v)\overline {A( v,w_{z,\zeta})}$ extends to an entire
function for each $z, \zeta$, and since the second factor is also
entire and does not vanish in a neighborhood of $(-z,-\zeta)\in
\C^d$, $A(u,v)$ extends also to an entire function. The fact that
it is of order $2$ follows from  (\ref{maj-g}), (\ref{below}) and
(\ref{above}).
\end{proof}
 We have completed the proof of Proposition  \ref{thBHamb}. \epf
\begin{proof}[Proof of Theorem \ref{thBHamb0}]
With the weaker assumption (\ref{HBA0}), we conclude  that
(\ref{prod-A})also holds, using the directional theorem for
Fourier transforms. We claim that $A(u,v)$ is an analytic function
of each of the variables $x$ and $y$. Indeed, as before, for every
function $w\in L^2(\R^d)$, the product  $A(u,v)\overline{A(v,w)} $
extends to a holomorphic function of $x$, $y$ being fixed, as well
as to a holomorphic function of $y$, $x$ being fixed. When
choosing $w$ as before, we conclude that the function extends to
an entire function of order $2$ in $x$, for fixed $y$. So, for
almost every fixed $y$, $A(u,v)(x,y)$ may be written as
$R_{y}(x)e^{2i\pi\scal{\alpha(y),x}}
e^{-\frac{\pi}{2}\norm{x}^2}$, with $R_{y}$ a polynomial of degree
at most $N$. It follows that the continuous function
$$e^{{\pi}(\norm{x}^2+\norm{y}^2)}|A(u,v)(x,y)|^2$$ is a
polynomial of degree at most $2N$ in each variable $x$, $y$. So it
is a polynomial, and the assumption (\ref{HBA}) is satisfied. We
can now use Proposition \ref{thBHamb} to conclude. \epf In
particular, it follows
 from Theorem \ref{thBHamb0} that there does not
exist two non zero functions $u,v$ which satisfy Condition
(\ref{HBA0}) for $N\leq d$. In  this case, where the conclusion is
that $u$ or $v$ is identically $0$, the condition can be relaxed
into a directional condition as for the case of the Fourier
transform.
\begin{corollary}
Assume that
\begin{equation}
\iint_{\R^{2d}}
\abs{A(u,v)}^2\frac {e^{\pi\abs{x_{j}}^2}}{(1+\abs{x_{j}})^M}
dxdy<+\infty
\hspace{2ex}\mathrm{and} \hspace{2ex}
\iint_{\R^{2d}}
\abs{A(u,v)}^2\frac {e^{\pi\abs{y_{j}}^2}}{(1+\abs{y_{j}})^N}
dxdy<+\infty
\label{HBA1}
\end{equation}
for some $j=1, \cdots,d$.  If $\min \{M,N\}=1$, then $u$ or $v$
vanishes. \label{corHBA}\end{corollary} \bpf When $d=1$, the
result follows from  Theorem \ref{thBHamb0}. Let us now consider
the case $d>1$. We assume that the condition is satisfied for
$j=1$. For $t'\in \R^{d-1}$ and for $f$ in $L^2(\R^d)$, we define
$f_{t'}(t_{1})=f(t_{1},t')$. It follows from Plancherel identity
that $$\int_{\R^{d-1}} \abs{A(u,v)(x,y)}^2dy'=\int _{\R^{d-1}}
\abs{A(u_{t'-\frac {x'}{2}},v_{t'+\frac
{x'}{2}})(x_{1},y_{1})}^2dy'\;.$$ Changing  variables, and using
Fubini's Theorem, it follows that, for almost every $\xi'$ and
$\eta'$ in $\R^{d-1}$, $$\iint_{\R\times\R}
\frac{\abs{A(u_{\eta'},v_{\xi'})(x_{1},y_{1})}^2}{(1+\abs{x_{1}}+\abs{y_{1}})^2}
e^{\pi(\abs{x_{1}}^2+\abs{y_{1}}^2)}dx_{1}dy_{1}<+\infty\;.$$ It
follows from the one dimensional case that either $u_{\eta'}$ or
$v_{\xi'}$ is identically $0$. So, for almost every $\xi'$ and
$\eta'$ in $\R^{d-1}$, $A(u_{\eta'},v_{\xi'})=0$. It follows that
$A(u,v)=0$. \epf
\medskip

As in section \ref{corollaires}, we may deduce from Theorem
\ref{thBHamb}\ a version of Hardy's theorem for the ambiguity
function. We state it here in a larger context, allowing general
positive definite quadratic forms. We do not give the proof, which
is obtained from a small modification of the proof of Theorem
\ref{thBHamb}. Let us remark that the constraints on degrees are
always elementary.

\begin{corollary} Let $u,v\in L^2(\R^d)$ and assume that
$$ \abs{A(u,v)(x,y)}\leq (1+\norm{x})^N(1+\norm{y})^Ne^{-\frac
12\scal{Cx,x}-\frac 12\scal{By,y}}. $$ We have the following
implications.
\begin{enumerate}
\item If $B-\pi^2 C^{-1}$ is positive, non zero, then either
$u$ or $v=0$\,.
\item If $B=\pi^2 C^{-1}$, there are polynomials $P,Q$ of
degree $\leq N$ and constants $\omega, a \in \R^d$ such that
$$
u(x)=P(x)e^{2i\pi\scal{\omega,x}}e^{-\scal{C(x-a),x-a}}
\mbox{ and }
v(x)=Q(x)e^{2i\pi\scal{\omega,x}}e^{-\scal{C(x-a),x-a}}.
$$
\end{enumerate}
\end{corollary}

\begin{remark} This corollary implies in particular that if
$u_0(x)=P_0(x)e^{-\alpha\norm{x}^2}$,  $v_0(x)=Q_0(x)e^{-\alpha\norm{x}^2}$ with
$P_0,Q_0$ polynomials and if $\abs{A(u,v)}=\abs{A(u_0,v_0)}$
then $u,v$ are of the form
\begin{align}
u(x)=&P(x)e^{2i\pi\scal{\omega,x}-\alpha\norm{x-a}^2}\notag\\
v(x)=&Q(x)e^{2i\pi\scal{\omega,x}-\alpha\norm{x-a}^2}\notag
\end{align}
with $P,Q$ polynomials and $\omega,a\in\R^d$.

The problem of finding $u,v$ from $u_0,v_0$ is known as the {\it radar
ambiguity problem} and has been considered by Bueckner \cite{Bu} and
de Buda \cite{dB} for $u_0,v_0$ as above. This remark
corrects the proof of \cite{dB}.

Further references on this problem may be found in \cite{Jaamb} and
\cite{GJP}.
\end{remark}

\medskip
Let us finally give a Gel'fand--Shilov type theorem.
\begin{theorem}
Let $1<p<2$, and let $q$ be the conjugate exponent. Assume that
$u,v\in L^2(\R^d)$ satisfy $$
\iint_{\R^d\times\R^d}|A(u,v)(x,y)|^2e^{2\pi\frac{a^p}{p}\abs{x_{j}}^p
}dxdy<+\infty \hspace{2ex}\mathrm{and} \hspace{2ex}
\iint_{\R^d\times\R^d}|A(u,v)(x,y)|^2e^{2\pi\frac{b^q}{q}\abs{y_{j}}^q}dxdy<+\infty
$$ for some $j=1, \cdots, d$ and for some positive constants $a$
and $b$. Then either $u$ or $v$ vanish if $ab> \abs
{\cos\left(\frac{p\pi}{2}\right)}^{\frac 1p}$. If  $ab< \abs
{\cos\left(\frac{p\pi}{2}\right)}^{\frac 1p}$, there exists a non
zero function $u$ such that the two conditions are satisfied by
$A(u,u)$. \label{GelShiAmb}
\end{theorem}
\bpf It is sufficient to consider the one dimensional case.
Otherwise, the proof works as in Corollary \ref{corHBA}. Let us
prove the first assertion. If we proceed as in  the proof of
Theorem \ref{thBHamb0}, we conclude at once that $$
A(u,v)(x,y)A(u,v)(-x,-y)=0 $$, using the similar result on Fourier
transforms. It remains to show that $A(u,v)$ is an analytic
function of each variable, which is obtained in the same way as
before using an auxiliary function $w$.

Let us now prove that for $ab< \abs {\cos\left(\frac{p\pi}{2}\right)}^{\frac 1p}$,
there exists a non
zero function $u$ such that the two conditions are satisfied by
$A(u,u)$. Using Plancherel formula, the first condition may as well be
written as
$$
\iint_{\R\times\R}|u\left(t-\frac{x}{2}\right)|^2|u\left(t+\frac{x}{2}\right)|^2e^{2\pi\frac{a^p}{p}\abs{x}^p
}dxdy<+\infty\;.
$$
Using the same change of variable as before, and the inequality
$$
| \eta-\xi|\leq 2^{p-1}(| \eta|^p+|\xi|^p)\,,
$$
we see that this integral is bounded by the square of the integral
$$
\int_{\R} |u(\xi)|^2 e^{4\pi\frac{(2^{1-2/p}
a)^p}{p}\abs{\xi}^p}d\xi\,.
$$
For the second integral, we use Lemma \ref{lem0} to write it in terms of
the Fourier transform of $u$. We obtain that it is bounded by
$$
\int_{\R} |\widehat u(\xi)|^2 e^{4\pi\frac{(2^{1-2/q}
a)^q}{q}\abs{\xi}^q}d\xi\,.
$$
The fact that there is a non zero function $u$ for which both integrals
are finite is an easy consequence of Theorem \ref{GelShi} and
Schwarz inequality, since  $\left(2^{1-\frac{2}{p}}a\right)
\times\left(2^{1-\frac{2}{q}}b\right)<
\abs {\cos\left(\frac{p\pi}{2}\right))}^{\frac {1}{p}}$.
\epf

\medskip
\begin{remark} At this stage, we would like to point out that we have not been
able to weaken the condition
$\iint\frac{\abs{A(u,v)(x,y)}^2}{(1+\norm{x}+\norm{y})^N}
e^{\pi(\norm{x}^2+\norm{y}^2)}dxdy<+\infty$ in Theorem
\ref{thBHamb0} into the condition $$
\iint\frac{\abs{A(u,v)(x,y)}^2}{(1+\norm{x}+\norm{y})^N}
e^{2\pi\abs{\scal{x,y}}}dxdy<+\infty \mbox{ or at least }
\iint\frac{\abs{A(u,v)(x,y)}^2}{(1+\norm{x}+\norm{y})^N}
e^{2\pi\norm{x}\norm{y}}dxdy<+\infty. $$ Also, we do not know
whether weaker conditions, with $\abs{A(u,v)(x,y)}$ in place of
its square,  are sufficient.
\end{remark}

\section*{Appendix}
We give here a simplified proof of Beurling-H\"ormander's Theorem,
which may be useful for elementary courses on Fourier Analysis. All
ideas are contained in Section 2.

We want to prove that a function
$f\in L^1(\R)$ which satisfies the inequality
\begin{equation}
    \iint_{\R\times\R}
|f(x)||\widehat{f}(\xi)|e^{2\pi\abs{x}\abs{\xi}}dxd\xi<+\infty
\label{condBHT1}
\end{equation}
is identically $0$. It is sufficient to show that the function
$g=e^{-\pi x^2}*f$ is identically $0$. Indeed, the Fourier
transform of $g$ is equal to $e^{-\pi \xi^2}\widehat{f}$. If it is
$0$, then $f$ vanishes also. Now $g$ extends to an entire function
of order $2$ in the complex plane. We claim that, moreover,
\begin{equation}
    \iint_{\R\times\R}
|g(x)||\widehat{g}(\xi)|e^{2\pi\abs{x}\abs{\xi}}dxd\xi<+\infty
\label{condBHT3}
\end{equation}
Indeed, replacing $g$ and $\widehat{g}$ by their values in terms
of $f$ and $\widehat{f}$ and using Fubini's theorem, we are led to
prove that the quantity $$\int_{\R}e^{-\pi [(x-y)^2
-2\abs{x}\abs{\xi}+2\abs{y}\abs{\xi} +\xi^2]} dx$$ is bounded
independently of $y$ and $\xi$. Taking $x-y$ as the variable, it
is sufficient to prove that $$\int_{\R}e^{-\pi [x^2
-2\abs{x}\abs{\xi} +\xi^2]} dx$$
is bounded by $2$, which follows
from the fact that $x^2-2\abs{x}\abs{\xi} +\xi^2$ is either $(x-\xi)^2$ or
$(x+\xi)^2$.

Now, for all $z\in \C$, we have the elementary inequality
$$\abs{g(z)}\leq
\int_{\R}|\widehat{g}(\xi)|e^{2\pi\abs{z}\abs{\xi}}d\xi\,,$$
so that there exists some constant $C$ such that
\begin{equation}
\int_{-\infty }^{+\infty}|g(x)| \sup_{\abs{z}=\abs{x}} |g(z)|dx\leq
C\,.
\label{condBHT4}
\end{equation}
We claim that the holomorphic function
$$G(z)=\int_{0}^z g(u)g(iu)du $$
is bounded by $C$. Once we know this, the end of the proof is immediate: $G$
is constant by Liouville's Theorem; so $g(u)g(iu)$ is identically
$0$, which implies that $g$ is identically $0$.

It is clear from (\ref{condBHT4}) that $G$ is bounded by $C$ on
the axes. Let us prove that it is bounded by $C$ for
$z=re^{i\theta}$ in the first quadrant. Assume that $\theta$ is in
the interval $(0,\pi/2)$. By continuity, it is sufficient to prove
that $$G_{\alpha}(z)=\int_{0}^z g(e^{-i\alpha} u)g(iu)du $$ is
bounded by $C$ for all $\alpha\in (0,\theta)$. But the function
$G_{\alpha}$ is an entire function of order $2$, which is bounded
by $C$ on the $y$-axis and on the half-line $\rho e^{i\alpha}$. By
Phagm\`en--Lindelh\"of principle, it is bounded by $C$ inside the
angular sector, which gives the required bound for
$|G_{\alpha}(z)|$. A similar proof gives the same bound in the
other quadrants.

\end{document}